\documentclass[psamsfonts]{conm-p-l}
\usepackage{amssymb}
\usepackage{pstricks}

\copyrightinfo{2008}{American Mathematical Society}

\newtheorem{theorem}{Theorem}[section]

\newtheorem{corollary}[theorem]{Corollary}

\theoremstyle{definition}
\newtheorem{definition}[theorem]{Definition}
\newtheorem{example}[theorem]{Example}

\theoremstyle{remark}
\newtheorem{remark}[theorem]{Remark}

\numberwithin{equation}{section}

\newcommand{\A}{{\mathcal A}}
\newcommand{\B}{{\mathcal B}}
\newcommand{\RR}{{\mathcal R}}
\newcommand{\VV}{{\mathcal V}}

\newcommand{\Z}{\mathbb{Z}}
\newcommand{\Q}{\mathbb{Q}}
\newcommand{\R}{\mathbb{R}}
\newcommand{\C}{\mathbb{C}}
\newcommand{\F}{\mathbb{F}}
\newcommand{\TT}{\mathbb{T}}
\newcommand{\p}{\mathbb{P}}
\newcommand{\K}{\mathbb{K}}
\newcommand{\G}{\Gamma}

\DeclareMathOperator{\Hom}{Hom}
\DeclareMathOperator{\Ker}{Ker}

\DeclareMathOperator{\gr}{gr}

\DeclareMathOperator{\codim}{codim}

\begin{document}

\title[Combinatorics and Topology of Hyperplane Arrangements]{%
A Survey of Combinatorial Aspects in the Topology of Complex Hyperplane Arrangements}

\author[A.~D.~Macinic]{Anca Daniela M\u acinic$^*$}
\address{Inst. of Math. Simion Stoilow\\
P.O. Box 1-764\\
RO-014700 Bucharest, Romania}
\email{Anca.Macinic@imar.ro}
\thanks{$^*$Partially supported by CNCSIS grant ID-1189/2009-2011 of the
Romanian Ministry of Education and Research.}


\begin{abstract}
We survey interactions between the topology and the combinatorics of complex hyperplane arrangements. 
Without claiming to be exhaustive, we examine in this setting combinatorial aspects of
fundamental groups, associated graded Lie algebras, higher homotopy groups, cohomology rings,
twisted homology with rank $1$ complex coefficients, and Milnor fibers.
\end{abstract}

\maketitle

\section{Introduction}
\label{sect:intro}

By a hyperplane arrangement $\A$ we understand in what follows either a collection $\{ H_1,\dots, H_n \}$ of codimension one subspaces
in a finite dimensional vector space $V$ or the subspace of the ambient space $V$ given by their union. We consider here only arrangements in {\em complex}
affine  spaces.

The starting point in the study of complex hyperplane arrangements was the braid 
arrangement, $\A=\{H_{ij}=\Ker(z_i-z_j)\;|\; 1\le i\neq j \le n \} \subset \C^n$. 
The results obtained for the complement $M_{\A}= \{z \in \C^n\;|\; z_i \neq z_j\; \forall i \neq j\}$ -- 
factorization of the Poincar\' e polynomial and the presentation of the cohomology algebra in terms 
of generators and relations  (Arnold \cite{Ar}), or the $K(\pi,1)$ property (Fadell-Neuwirth \cite{FN}) --
initiated typical directions of investigation in arrangement theory.
Resolving a conjecture of Arnold, Brieskorn showed that the cohomology ring of the complement,
$M_{\A}= V \setminus \bigcup_{i=1}^n H_i$, of an
arbitrary complex arrangement $\A$  is generated by the classes of the closed $1$-forms 
$\frac{1}{2 \pi \sqrt{-1}}(\frac{d \alpha_H}{\alpha_H})$, where $\alpha_H=0$ are the defining equations for
the hyperplanes $H \in \A$.

An emblematic result for the field is the theorem of Orlik and Solomon \cite{OS}, that proves 
an isomorphism 
between the cohomology ring of the complement, $H^*(M_\A)$, and  the {\it Orlik-Solomon algebra} 
of the arrangement, $A^*(\A)$, which is  determined by the lattice 
$\mathcal L(\A)=\{\bigcap_{H \in{\B}}H\;|\;\B \subset \A\}$, ordered by reverse inclusion.

We say that a topological invariant (property) of the complement $M_{\A}$ is {\em combinatorial} 
if it depends only on the intersection lattice $\mathcal L(\A)$.
  
We have chosen to survey here several natural, interesting topological facets of the 
combinatorial determination problem (based on personal taste, and with no intention of
exhausting the subject).
Furthermore, ample existing expository work covers various aspects of the field of hyperplane arrangements, see for instance \cite{FR1}, \cite{FR2}, \cite{Fa0}, \cite{S1}, \cite{Y1}.

A brief discussion of our main topics follows.

Projection onto the first coordinates, $\C^n \to \C^{n-1}$, restricts to a topological fibration
between the corresponding braid complements, with fiber $\C \setminus \{ n-1 \; \text{points} \}$.
Arrangements whose complement may be obtained from the simplest arrangement, $\{ 0\}\subset \C$,
by iterating a similar fibration construction, were introduced by Falk-Randell \cite{FR0} and called
{\em fiber-type}. A remarkable result (due to Terao \cite{Te} and described in Section \ref{sect:lattice})
says that this is equivalent to $\mathcal L(\A)$ being {\em supersolvable}, in the sense of Stanley \cite{St}.

In Section \ref{sect:cohomology}, we examine the converse of the Orlik-Solomon theorem: the
cohomology ring $H^*(M_\A)$, together with its {\em natural $1$-marking} (in the sense of
Choudary-Dimca-Papadima \cite{CDP}), determines the lattice $\mathcal L(\A)$.

Section \ref{sect:pi1} treats the homotopy groups, $\{ \pi_i(M_{\A})\}_{i\ge 1}$. As shown by
Rybnikov \cite{Ryb}, the {\em fundamental group} $\pi_1(M_{\A})$ is {\em not} combinatorial (see also \cite{AAB0}, \cite{AAB1}).
Nevertheless, its rational {\em associated graded Lie algebra}, $\gr_*(\pi_1(M_{\A}))\otimes \Q$,
is combinatorial, by a result of Kohno \cite{Ko}. The nilpotence of $\pi_1(M_{\A})$ is
equivalent to the property of $\mathcal L(\A)$ of  being in general position
in codimension $2$; this fact is the first step in the proof of a famous conjecture of Zariski.

The $K(\pi,1)$ property of $M_{\A}$ means by definition the vanishing of the higher homotopy
groups $\pi_i(M_{\A})$, for $i>1$. It is easy to see that fiber-type implies  $K(\pi,1)$, but the
converse is not true, as shown by an example of Falk \cite{Fa1}. A basic result of Hattori
\cite{Ha} implies that the $K(\pi,1)$ property does not hold, if the lattice $\mathcal L(\A)$ is
{\em generic}. Jambu-Papadima introduced in \cite{JP1} the (combinatorially defined)
{\em hypersolvable} class of arrangements, which (strictly) contains both the fiber-type 
and the generic arrangements. They also proved that the $K(\pi,1)$ property is combinatorial,
within the hypersolvable class. In general, this question is open.

Section \ref{sect:special} is devoted to the hypersolvable class. For the {\em split solvable} subclass,
we highlight a remarkable perfect match, namely the fact that the lattice $\mathcal L(\A)$ and 
the homeomorphism type of the complement $M_{\A}$ determine each other \cite{CDP}.

Rank $1$ complex locally constant sheaves (alias {\em local systems}) on  $M_{\A}$ may also be
viewed as characters $\rho$ belonging to the character torus, $\TT_{\pi_1(M_{\A})}=(\C^*)^n$.
Combinatorial determinacy of the {\em twisted cohomology} $H^*(M_{\A}, {}_{\rho}\C)$ is
a major open problem in the field. We review in Section \ref{sect:VR} 
several known facts in this direction. 

The {\em Milnor fiber} of a degree $n$ homogeneous polynomial $f\in \C [x_1, \dots, x_l]$ is
an important object in singularity theory, introduced by Milnor in \cite{Mi}. It is defined by
$F_f := f^{-1}(1)$, and comes endowed with a natural automorphism, induced by multiplication in
$\C^l$ by a primitive $n$-root of unity. The induced action on $H_*(F_f, \Q)$ is called the
{\em algebraic monodromy}. Its study is a very active area of research in singularity theory.
In arrangement theory, the Milnor fiber of $\A$, $F_{\A}$, is associated to the defining polynomial of
$\A$, $Q(\A):= \Pi_{H\in \A} \alpha_H$.

A long-standing open problem in this context is the combinatorial determination of the algebraic 
monodromy action on $H_1(F_{\A}, \Q)$. We address this question in Section \ref{sect:Milnor}.
Following \cite{MP}, we present an affirmative answer, in terms of combinatorial objects 
defined in positive characteristic, for subarrangements of type $B$ Coxeter arrangements,
parametrized by {\em graphs}. We also emphasize a general relation between algebraic monodromy
and twisted cohomology, which is key to the proof.

\section{Lattice and complement of an arrangement}
\label{sect:lattice}

Assume $\A \subset V,\;V$ a complex vector space of dimension $l$, is an arrangement of hyperplanes.
Choosing a system of coordinates $x_1, \dots, x_l$ we may  express each hyperplane 
$H \in \A$ as the zero set of a degree $1$ polynomial $\alpha_H \in \C[x_1, \dots, x_l]$. 
The product $Q(\A)=\Pi_{H \in \A}\alpha_H$ is called the defining polynomial of $\A$. The
{\em complement} of $\A$ is $M_{\A} := \C^l \setminus Q(\A)^{-1}(0)$. 

 For arrangements with the property that $\bigcap_{H \in \A}H \neq \emptyset$, by a choice of 
 coordinates, one can assume that all hyperplanes contain the origin, hence all 
 corresponding equations $\alpha_H$ are homogeneous. Such an arrangement is called {\it central}.
 
Recall that the {\em intersection lattice} $\mathcal L(\A)$, i.e., the set of all subspaces 
which are intersections of hyperplanes 
of $\A$, is a poset with a rank function, whose partial order is
given by reverse inclusion, and rank given by $r(\cap_{H\in \B} H) := \codim (\cap_{H\in \B} H)$. 
The minimal element of $\mathcal L(\A)$ is $V$, seen as intersection 
of the empty set of hyperplanes.
\begin{example}
\label{latice}
The lattice of the {\it boolean} arrangement $\A \subset \C^l,\; Q(\A)=\Pi_{i=1}^l x_i$, contains 
$2^l$ subspaces $H_{i_1} \cap \dots \cap H_{i_s}$, indexed by all ordered subsets 
$i_1 < i_2< \dots <i_s$ of $\{1, 2, \dots, l\}$, with rank $s$. In this case, $M_{\A}= (\C^*)^l$.
\end{example}
 
A basic combinatorial invariant of an arrangement is the {\it rank}, denoted $r(\A)$, defined as 
the codimension of a maximal element of $\mathcal L(\A)$. It is  well defined, as we can see from \cite[Lemma 2.4]{OT}.  
The levels $\mathcal L_p(\A), \;0 \leq p \leq r(\A)$ consist of the elements of rank $p$
of the lattice $\mathcal L(\A)$. $\A$ is called {\it essential} if $r(\A)=l$. We can always assume 
that a central arrangement
is essential, without changing the homotopy type of $M_{\A}$, eventually seeing it 
as an arrangement in the quotient of the ambient
vector space $V$ by the {\em center} of the arrangement, $\cap_{H\in \A} H$.

\begin{example}
\label{essential}
\begin{enumerate}
\item \label{boolean}
The {\em boolean} arrangement is an essential arrangement.
\item \label{converse}
The {\em braid} arrangement $\A \subset \C^l,\; Q(\A)=\Pi_{1 \leq i<j \leq l}(x_i-x_j)$ 
has rank $l-1$, 
hence is not an essential arrangement.
\end{enumerate}
\end{example}

\begin{definition}
\label{cone}
The {\em cone} over an arbitrary arrangement $\A \subset \C^l$ is the  
central arrangement $c\A \subset \C^{l+1}$ 
with the defining polynomial $Q(c\A)=x_0 \overline{Q(\A)}$, where $\overline{Q(\A)} \in \C[x_0 \dots x_l]$
is the homogeneization of $Q(\A) \in \C[x_1 \dots x_l]$.
Conversely, we can apply a deconing procedure to a central arrangement $\A \subset \C^{l+1}$ to obtain 
an affine arrangement $d\A \subset \C^l$. Choose coordinates such that $\A$ contains a hyperplane $H$
of equation $x_0=0$. Then $Q(d\A)$ is the evaluation of $Q(\A)$ at $(1, x_1, \dots, x_l)$.
\end{definition}

There is a topological reason behind this construction, emerging from projectivisation. 
The Hopf fibration $p: \C^{l+1} \setminus \{0\} \longrightarrow  \C \p^{l}$, with fiber $\C^*$, sends 
the complement $M_{c\A}$ to $M_{\A}$ and produces a decomposition $M_{c\A} \cong M_{\A} \times \C^*$,
or, equivalently,  for central $\A$, $M_{\A} \cong M_{d\A} \times \C^*$(see \cite[Proposition 5.1]{OT}). 
Here, we  see the complement of an arbitrary arrangement $\A \subset \C^l$ as the complement of 
the projectivised arrangement $\p (c\A) \subset \C\p^l$: first remove the hyperplane at infinity, 
$H_0= \{ x_0=0\}$, 
then the rest of the hyperplanes from $\C\p^l \setminus H_0 \cong \C^l$.

\begin{definition}
\label{modular}
Let $\A$ be central. An element $X \in \mathcal L(\A)$ is called {\em modular} if 
$X+Y \in \mathcal L(\A)$, for any $Y\in \mathcal L(\A)$.
For other equivalent definitions see \cite[\S2.1]{OT}.
\end{definition}

\begin{example}
\label{modular1}
All hyperplanes are modular elements. The center of the lattice is also modular. 
All elements in the lattice of the boolean arrangement are modular.
\end{example}

\begin{definition}[Stanley \cite{St}]
\label{stanley}
A central arrangement $\A\subset \mathbb{C}^{l}=V$  is called {\em supersolvable} if there is a 
maximal chain of modular elements:
\begin{equation}
\label{eq:21}
V=X_{0}\supset X_{1}\supset\dots \supset X_{r}=C(\A)
\end{equation}
where $C(\A)= \bigcap_{H \in \A}H$ is the center.
\end{definition}

\begin{example}
\begin{enumerate}
\item \label{bool3}
The boolean arrangement from Example \ref{latice} is supersolvable, 
with a maximal chain of modular elements given for instance by 
$V \supset \{x_1=0\}\supset \{x_1=x_2=0\} \dots \supset \{x_1= \dots=x_l=0\}$.
\item \label{braidchain}
A chain of modular elements of length $l-1$ for the braid arrangement in Example 
\ref{essential} (\ref{converse}) is the following: 
$V \supset \{x_1=x_2\} \dots \supset \{x_1= \dots=x_l\}$.
\end{enumerate}
\end{example}

The next definition gives a topological interpretation 
of supersolvability: as shown by Terao \cite{Te}, $\mathcal L(\A)$ is supersolvable
if and only if $\A$ is fiber-type.

\begin{definition}[Falk-Randell \cite{FR0}]
\label{fibertype}
An essential arrangement $\A \subset \mathbb{C}^{l}$  is said to be {\em fiber-type} if there is 
a sequence of subarrangements
in $\A$ of length $l$ :
\begin{equation}
\label{eq:22}
\{ H \}=\A_{1} \subset \A_{2} \subset \dots \subset \A_{l}=\A
\end{equation}
so that there is a tower of topological fibrations $M_{\A_{i}}\longrightarrow M_{\A_{i-1}}$, 
induced by linear maps, with fiber $F_i$ a complex line with a finite number of points removed.
\end{definition}

Fiber-type arrangements are a class of arrangements very present in the literature,
for which significant topological information is available.

\section{Cohomology algebra and twisted cohomology}
\label{sect:cohomology}

Let $X$ be an arbitrary topological space and $\K$ a unitary commutative ring. 
The  {\em cohomology algebra} of $X$ with (untwisted) coefficients in $\K$ is 
the graded ring 
$H^*(X, \K)=\oplus_{i \geq 0}H^i(X, \K)$ with a multiplication 
$H^i \otimes H^j \overset{\cup} \rightarrow H^{i+j}$, called cup product,
which is commutative in the graded sense, that is, $a \cup b= (-1)^{ij} b\cup a$,
for $a\in H^i$ and $b\in H^j$. See for instance \cite{H}.
By convention, $H^*(X):=H^*(X,\Z)$.

\begin{example}
\label{deRham}
The cohomology of a differentiable ($\mathcal{C}^{\infty}$) manifold $X$ with real coefficients, 
$H^*(X,\R)$, is isomorphic to the cohomology of the de Rham algebra $\Omega^*_{dR}X$, where 
the cup product is induced by the exterior product on differential forms.
\end{example}

\begin{remark}
\label{lineardecomp}
When $X=M_{\A}$ is the complement of a fiber-type arrangement, we have an additive combinatorial 
decomposition of $H^*(M_{\A})$. More precisely, in the notations of Definition \ref{fibertype}:
\begin{equation}
\label{factor}
H^*(M_{\A}) \cong H^*(F_1) \otimes \dots \otimes H^*(F_l)\, ,
\end{equation}
as graded abelian groups.
\end{remark}

Naturally, one may ask whether combinatorial determination holds for the cohomology algebra of 
an arbitrary arrangement complement. The answer is positive, as we shall see in Theorem \ref{thm:OS}. 
To state the result, we need an essential combinatorial ingredient, which we describe next.

Denote by $E^*(\A)$ the exterior algebra generated, over an arbitrary unitary commutative ring 
$\mathbb{K}$ by elements $e_H$ of degree $1$, one for each $H \in \A$. The 
{\it Orlik-Solomon algebra} of $\A,\; A_{\mathbb{K}}^*(\A)$ is defined as a quotient of 
 $E^*(\A)$ by an ideal $I$ generated by dependency relations in $\mathcal L(\A)$. More precisely,
 one has a canonical boundary operator, $\partial \colon E^*(\A)\to E^{*-1}(\A)$, defined by
 \begin{equation}
 \label{eq:32}
 \partial (e_{H_{i_1}}\cdots e_{H_{i_t}})= \sum_{j=1}^t (-1)^{j-1} e_{H_{i_1}}\cdots 
 \widehat{e_{H_{i_j}}} \cdots e_{H_{i_t}}\, ,
 \end{equation}
 where $\widehat{\bullet}$ denotes deletion of $\bullet$.
 For central arrangements, $I$ is generated by boundaries of elements $e_{H_{i_1}} \dots e_{H_{i_t}}$ 
 such that $r(\bigcap_{j=1}^{t}H_{i_{j}}) <t$. For arbitrary affine arrangements, one
 also assumes in the above construction that $\bigcap_{j=1}^{t}H_{i_{j}}\ne \emptyset$, and then
 also adds to the generating set of $I$ elements $e_{H_{i_1}} \dots e_{H_{i_t}}$ such that 
 $ \bigcap_{j=1}^{t}H_{i_{j}}= \emptyset$. We denote the generators of $A_{\mathbb{K}}^*(\A)$ 
 again by $e_{H}$, for  $H \in \A$.
 
 The simplest example of OS-algebra is the one associated to the boolean arrangement. Since 
 there are no dependency relations among the hyperplanes,  the OS-algebra 
 is just the exterior algebra, $E^*(\A)$.
 
 \begin{theorem}[Orlik-Solomon \cite{OS}]
 \label{thm:OS}
 There is an isomorphism of graded algebras
 \begin{equation}
 A_{\mathbb{K}}^*(\A) \cong H^*(M_{\A}, \mathbb{K}),
 \end{equation}
 that is, $\mathcal L(\A)$ determines $H^*(M_{\A}, \mathbb{K})$ as a graded ring. 
\end{theorem}
 
 A similar object is the {\it quadratic Orlik-Solomon algebra} of a central arrangement 
 $\A$,  denoted $\overline{A^*_{\K}}(\A)$, defined as the  quotient of $E^*(\A)$ 
 by the ideal spanned by relations of type
\begin{center}
 $e_{H_{i}}e_{H_{j}}+e_{H_{j}}e_{H_{k}}=e_{H_{i}}e_{H_{k}}\, ,$
\end{center}
for $H_{i},H_{j},H_{k}$ such that $r(\{H_{i},H_{j},H_{k}\})=2$.
For  fiber-type arrangements, the Orlik-Solomon algebra 
is {\it quadratic}, i.e., $\overline{A^*_{\K}}(\A)=A^*_{\K}(\A)$; see Falk \cite{Fa1}
and Shelton-Yuzvinsky \cite{SY}.

  The converse of Theorem \ref{thm:OS} is however no longer true:
  the arrangements from Example 3.77 \cite{OT} have isomorphic cohomology algebras, but
  different lattices. This inconvenient disappears, when adding information provided by
  {\it natural $1$-markings}. 
  
\begin{definition}
\label{1mrk}
A {\em $1$-marking} of a space $X$ of finite type is a distinguished ordered $\Z$-basis of $H^1(X)$.
A {\em $1$-marking} of a finitely generated group $G$ with free abelianization 
$G_{ab}$ is a distinguished ordered 
$\Z$-basis of $G_{ab}$. {\em $1$-marked isomorphisms} must preserve the 
distinguished ordered $\Z$-bases, in both cases.
\end{definition}

Ordered hyperplane complements have a {\em natural} $1$-marking in cohomology, namely
$\{ e_H \mid H\in \A \}$.

\begin{theorem}[Kawahara \cite{K}, Choudary-Dimca-Papadima \cite{CDP}]
\label{thm:1m}
For central arrangements, the cohomology ring $H^*(M_{\A})$, together with its natural $1$-marking,
determines the lattice $\mathcal L(\A)$.
\end{theorem}

Now, let us recall another  important topological invariant of a path-connected space $X$. 
Take the set of loops on $X$, starting and ending at some arbitrary 
fixed point, and identify  two loops if one can be continuously deformed into another. This set,
together with the obvious multiplication given by concatenation of loops, becomes a group,
the {\it fundamental group} of the space $X$, denoted by $\pi_1(X)$.

\begin{example}
\label{wedgeS1}
The fundamental group of the bouquet of $n$ circles, $\bigvee_{n}S^1$, is 
the free group on $n$ generators, $\F_n$.
\end{example}

 If $\A$ is an ordered arrangement, it has a {\em natural} $1$-marking 
 of $\pi_1(M_{\A})$, since $(\pi_1(M_{\A}))_{ab}$ is the $\Z$-dual of
 $ H^1(M_{\A})$. 
 
 Let $M$ be a finite type connected CW-complex with fundamental group $G$. By a {\it local system} 
 on $M$ we mean a right $\Z[G]$-module, over the group ring $\Z[G]$.
 
\begin{example}
\label{rk1}
The {\em character torus} $\TT_{G}:=\Hom (G, \C^*)$ is an algebraic variety that encodes the
so-called {\em rank $1$ complex local systems} on $M$. 
We shall denote by $_{\rho}\C$ the $\Z[G]$-module structure on $\C$ associated to an element
$\rho \in \TT_{G}$. 
\end{example}

\begin{example}
\label{ex:exp}
When $M=M_{\A}$ is the complement of an arrangement $\A$ with $n$ hyperplanes, one has 
that $\TT_{G}=(\C^*)^n$, and the exponential map $\C \overset {\exp} \longrightarrow \C^*$ 
induces an analytic morphism 
$\C^n= H^1(M_{\A},\C) \overset {\exp} \longrightarrow  \TT_{G}=(\C^*)^n$.

The  local system associated to a character $\rho \in \TT_{G}=\Hom (G_{ab},\C^*)$ is called 
{\it equimonodromical} if it is constant on the natural basis of $G_{ab}$.

 A {\it system of complex weights} 
 $\alpha=(\alpha_1, \dots, \alpha_n), \; \alpha_i \in \C,\; 1 \leq i \leq n$, on the hyperplanes 
 of $\A$, or equivalently an element $\alpha_{1}e_1+ \dots \alpha_ne_n \in A^1_{\C}(\A)$, 
 exponentiates to an element of the character torus, $\rho=(\exp(\alpha_1), \dots, \exp(\alpha_n))$.
 If $\alpha_i=2 \pi \sqrt{-1}\frac{m_i}{d},\; 1 \leq i \leq n$, for some integers $m_1, \dots, m_n$ 
 with $gcd(m_1, \dots m_n)=1$ and $d \in \Z_{>0}$, the character is called {\em rational}.
 Rational equimonodromical characters are very important in the study of 
 the Milnor fiber cohomology (see Section \ref{sect:Milnor}).
\end{example}

Let $(\tilde{C}_*, \tilde{d}_*)$ denote the equivariant chain complex of the universal cover of $M$.
It is a free chain complex of finitely generated left $\Z[G]$-modules.
The {\em twisted}  homology (respectively cohomology) with respect to a local system $L$ 
is the homology of $L\otimes_{\Z[G]} (\tilde{C}_*, \tilde{d}_*)$
(respectively the cohomology of $\Hom_{\Z[G]} ((\tilde{C}_*, \tilde{d}_*), L)$). For cohomology,
we  turn the right $\Z[G]$-module $L$ into a
left $\Z[G]$-module, in the usual way. Notation: $H_*(M, L)$ (respectively $H^*(M, L)$).

For example, the trivial character $\rho=(1, \dots, 1)$ gives rise to the usual cohomology with 
(untwisted) complex coefficients, that is, $H^*(M_{\A}, \;_{\rho}\C) = H^*(M_{\A}, \C)$, for $\rho =1$.

\begin{example}
\label{points}
Let $\A$ be an arrangement of $n$ points in $\C$. The complement has the homotopy type of a 
bouquet of $n$ circles, a $K(\pi,1)$ space (see Definition \ref{high}), with fundamental group 
the free group on $n$ generators. Then $H^*(M_{\A}, \;_{\rho}\C)=0$, for $n=1$ and $\rho \ne 1$.
\end{example}

{\it Having Theorem \textup{\ref{thm:OS}} in mind, we may ask: is there an analog 
of the Orlik-Solomon algebra for twisted coefficients, in other words, is the twisted 
cohomology combinatorially determined, at least for rank $1$ complex local systems?
This appears to be the major open problem in the topology of arrangements.} See
Section \ref{sect:VR} for more on this.
 
\section{Homotopy groups}
\label{sect:pi1}

We discuss first the fundamental group of 
the complement of an arrangement $\A$. As before, we are interested in combinatorial aspects
of $\pi_1(M_{\A})$.

For instance, the abelianity of the fundamental group of the complement is equivalent to 
its nilpotency, and equivalent to the property of the lattice of being 
in general position in codimension $2$. See \cite{Mac} for more equivalences.

\begin{example}
\label{fiber:pi1}
\begin{enumerate}
\item \label{Ftype}
The fundamental group of the complement of a fiber-type arrangement is 
an iterated almost-direct 
product of free groups on a finite number of generators (those free groups are exactly 
the fundamental groups of the  fibers $F_i$ from Definition \ref{fibertype}).
\item \label{Splitsolv}
If $\A$ is  split solvable of type $m=(m_1, \dots, m_r)$ (see Section \ref{sect:special}
for the definition), then $\pi_1(M_{\A})$ is isomorphic to the product of free groups 
$\F_{m_1} \times \dots \times \F_{m_r}$.
\end{enumerate}
\end{example}

\begin{remark}
\label{lefschetz}
Concerning the fundamental group, notice that line arrangements in $\C\p^2$ are the general case, 
via the Lefschetz hyperplane section theorem. See \cite{DP} for more results in
this direction.
\end{remark}

Although algorithms that give a presentation for $\pi_1(M_{\A})$ have been developed, 
it is well known that $\pi_1(M_{\A})$ is not a combinatorial invariant (Rybnikov, \cite{Ryb}). 
Conversely, one may ask: does $\pi_1$ determine the intersection lattice?
This time the answer is affirmative, for line arrangements, as long as we deal with
the naturally $1$-marked  fundamental group of the complement. 
More generally, for an arbitrary projective arrangement $\A$, the lattice is recovered 
up to some level (at least $2$) from $\pi_1(M_{\A})$ endowed with the natural
$1$-marking (see \cite[Theorem 2.1]{CDP}).

 Let $G$ be a group. Define a descending series of normal subgroups of $G$, called the
 {\it lower central series} associated to the group $G$, by $G_1:=G$ and 
 $G_i:=[G_{i-1},G]$ for $i \geq 2$, where $[x,y]:= xyx^{-1}y^{-1}$ is the group commutator.
 Then, $\gr_*(G) =\bigoplus_{i \geq 1} (\gr_i(G):= \frac{G_i}{G_{i+1}})$ has a natural 
 graded Lie algebra structure with the Lie bracket induced by the commutator. The
 {\it rational associated graded Lie algebra} is  the graded 
 Lie algebra over $\Q,\; \gr(G) \otimes \Q$.

\begin{example}
\label{lie}
\begin{enumerate}
\item \label{Zn}
If $G=\Z^n$ is the free abelian  group on $n$ generators, then $\gr_*(G)=\gr_1(G)$ is
abelian, i.e., $[,]=0$.
\item \label{Fn}
If $G=\F_n$ is the free group on $n$ generators, then $\gr(G)$ is the
free Lie algebra on $n$ generators.
\end{enumerate}
\end{example}

In the arrangement case, we point out that the rational 
associated graded Lie algebra of $\pi_1(M_{\A})$ is combinatorial, 
as a consequence of the formality of the complement space $M_{\A}$. See Kohno \cite{Ko}.
{\it It is an open question whether $\gr(\pi_1(M_{\A}))$ is determined by the lattice $\mathcal L(\A)$,
since it is known to contain  torsion elements.}

Next, we recall the definition of other significant objects in algebraic topology, 
that play an important role in the homotopy classification of topological spaces:
the {\em higher homotopy groups}. The $n$-th homotopy group of
a topological space $X, \; \pi_n(X)$, is the set of continous maps $[0,1]^n \rightarrow X$ 
that take the boundary of the $n$ dimensional cube $[0,1]^n$ to a fixed point in $X$, 
the base point, where two such maps are identified if one 
can be continuously deformed into another, by a deformation that fixes the base point. 
To describe the group operation by analogy with the case of the first homotopy group 
(the fundamental group), note that "the concatenation" of two $n$-cubes in $X$ means 
that they are glued along a face. For $n\ge 2$, $\pi_n(X)$ is abelian, 
yet very hard to compute.

\begin{definition}
\label{high}
A topological space $X$ is called $K(\pi,1)$ (or {\em aspherical}) if all homotopy groups
$\pi_{\geq 2}(X)$ are trivial.
\end{definition}

\begin{example}
\begin{enumerate}
\item \label{wedgeS12}
The bouquet of $n\geq 1$ circles, $\bigvee_{n}S^1$, is a $K(\pi,1)$ space.
\item \label{Kpi2}
The complement of a fiber-type arrangement is also aspherical.
\end{enumerate}
\end{example}

 Deligne showed that a large family of arrangements, the {\it complexified  simplicial arrangements}
 (\cite{De}) have $K(\pi,1)$ complements. Among them are the Coxeter arrangements, 
 which we describe next.
 A {\it Coxeter group} is a finite group generated by real, orthogonal reflections. 
 Irreducible Coxeter groups are classified. Their list is: 
 $A_{l},\;B_{l},\;D_{l},\;E_{6-8},\;F_{4},\;G_{2},\;H_{3-4},\;I_{2}(p)$. A {\it Coxeter arrangement}
 is the set of reflecting hyperplanes of the reflections in a Coxeter group.

Defining polynomials for Coxeter arrangements of type $A_l,  B_l$ and $D_l$ are, in this order:
$\Pi_{1\le i <j\le l}(x_i-x_j)$, $\Pi_{1\le i< j\le l}(x_i\pm x_j)\Pi_{i \in \overline{1,l}}x_i$, 
respectively $\Pi_{1\le i <j\le l}(x_i \pm x_j)$.

Examples of arrangements with non-$K(\pi,1)$ complements are also known.

\begin{definition}
\label{counterexp}
Let $\A \subset \C^r, \;r>2$, be a central arrangement with $n>r$ hyperplanes. 
$\A$ is called {\em generic} if any subarrangement $\mathcal{B} \subset \A$ 
of $r$ hyperplanes has rank $r$. 
\end{definition}

 A fundamental result, due to Hattori \cite{Ha}, on the homotopy type of the complement of a generic 
 arrangement, implies $\pi_{r-1}(M_{\A}) \neq 0$, in other words generic arrangements 
 are never $K(\pi,1)$, for $r \geq 3$.
  
 {\it So, one may ask: is the $K(\pi,1)$ property combinatorial? This is another interesting
  open question.}
  The answer is positive for the class of hypersolvable arrangements
   (defined in Section \ref{sect:special}). 
   
\begin{theorem}[Jambu-Papadima \cite{JP1}]
\label{test}
Let $\A$ be a hypersolvable arrangement. Then $M_{\A}$ is $K(\pi,1)$ if 
and only if $\mathcal L(\A)$ is supersolvable.
\end{theorem}

 The hypersolvable class includes both generic and supersolvable arrangements, and more. 
 For instance, it also contains the
 split solvable arrangements (also defined in the next section), 
 which are in general neither generic nor supersolvable. 

\section{Hypersolvable arrangements}
\label{sect:special}

Jambu-Papadima \cite{JP1} introduced hypersolvable arrangements. This  class  generalizes  
the fiber-type class, preserving many of its topological and combinatorial properties.

Inspired by the combinatorics of the braid arrangement, the conditions in the inductive definition 
of hypersolvability involve only collinearity relations, or equivalently elements in the lattice 
of the arrangement up to rank $2$.

The building blocks are certain pairs of central arrangements,
$(\A, \B),\; \B \subset_{\neq} \A$, 
with the convention that elements in $\B$ are denoted by $\alpha, \beta, \gamma, \dots$ 
and elements in $\overline{\B}:=\A \setminus\B$ by $a, b, c, \dots$. We will identify the hyperplanes
in $\A \subset V$ with points in $\mathbb{P}(V^*)$.

\begin{definition}[Jambu-Papadima \cite{JP1}]
\label{exthyper}
An arrangement pair 
$(\A, \B)$ is a {\em solvable extension} if the following conditions are satisfied: 
\begin{enumerate}
\item \textbf{closure}: any distinct points $\alpha, \beta, a$, with
$\alpha, \beta \in \mathcal{B}$ and $a \in \overline{\mathcal{B}}$ are non-collinear;
\item \textbf{completion} : for any two distinct points $a, b \in \overline{\mathcal{B}}$, 
there is a point $\alpha \in {\mathcal{B}}$ on the line determined by $a$ and $b \,$;
if the closure condition holds, $\alpha=:f(a,b)$ is uniquely determined;
\item \textbf{solvability}: for any distinct $a,b,c$ in $\overline{\mathcal{B}}$, 
the three points $f(a,b),f(c,b)$,
$f(a,c)$ are either equal or collinear.
\end{enumerate}
\end{definition}

\begin{definition}
\label{hyper}
A central arrangement $\A$ is called {\em hypersolvable} if it admits a 
composition series of solvable extensions:
$$
\{H\}=\A_{1}\subset \A_{2}\subset \cdots \A_{i}\subset \A_{i+1} \cdots \subset\A_{l}=\A
$$
\end{definition}
  For example, for generic arrangements, known to be hypersolvable, a composition series 
  can be easily constructed starting with an arrangement containing just one hyperplane, and
  adding one more hyperplane for each extension.
  
As in the fiber-type case, the fundamental group $\pi_1(M_{\A})$, for $\A$ hypersolvable, is 
an iterated almost direct product of finitely generated free groups 
$\mathbb{F}_{m_i}, \; m_i=|\A_i \setminus \A_{i-1}|$. 

 The length $l=l(\A)$ of a composition series, called {\it the length of $\A$}, 
 is a combinatorial invariant:
 it coincides with the degree of the Poincar\' e polynomial of 
 the quadratic OS-algebra $\overline{A^*}(\A)$, 
 see \cite[Theorem B]{JP1}.
 
\begin{theorem}[Theorem D, \cite{JP1}]
\label{thm:lr}
Let $\A$ be a hypersolvable arrangement. Then $\A$ is supersolvable 
if and only if $l(\A)=r(\A)$.
\end{theorem}

 Hence comparing the length of a hypersolvable arrangement to its rank is a simple
 combinatorial test for the asphericity of the complement.

 As a general rule, $l(\A) \geq r(\A)$ for $\A \subset V$ hypersolvable. Each solvable extension 
 $\B \subset_{\neq} \A$ increases the rank with at most $1$, that is $r(\A) \leq r(\B)+1$. 
 Accordingly, singular ($r(\A) = r(\B)$) or nonsingular
 ($r(\A) =r(\B)+1$) extensions appear. Singular extensions may be deformed, by 
 constructing an arrangement $\tilde{\A} \subset V \times \C,\; r(\tilde{\A})= r(\B)+1$, with 
 the same fundamental group of the complement as $\A$ and the same 
 collinearity relations.
 In the end this process produces for an arbitrary hypersolvable $\A$ a fiber-type 
 deformation $\tilde{\A} \subset V \times \C^s$, where $s$ is the number
 of singular extensions in a composition series of $\A$ (\cite{JP2}).
 
Another large family of examples of hypersolvability is described next.
 {\it Split solvable} arrangements are line arrangements $\A=\{H_0, \dots,H_n\} \subset \C\p^2$ 
 with a simple intuitive combinatorics. Take
 $H_0$ to be the line at infinity.  Let $p_1, \dots, p_r$ be the intersection points on $H_0$,
 with multiplicities $m_1+1, \dots, m_r+1$. By definition, $\A$ is split solvable if
 it has only double intersection points outside $H_0$.
  Equivalently, the affine arrangement in $\C\p^2\setminus H_0 \cong \C^2$ can be pictured 
  as $n$ lines with $r$ parallel directions, each direction containing $m_i$ lines, 
  $1 \leq i \leq r$, and $d\A$ has only double intersection points.
  
Call $m=(m_1, \dots, m_r)$ the {\it combinatorial type} of $\A$. Obviously the type $m$ 
completely describes the combinatorics of the arrangement.
A concrete type $m$ split solvable arrangement is constructed in \cite{CDP} in the following way.
Consider the arrangement $\hat{\A}(m)$ in $\C\p^r$ of equation
\begin{equation}
\label{longeq}
 z_0(z_1-z_0)(z_1-2z_0) \dots (z_1-m_1z_0)(z_2-z_0)\dots (z_2-m_2z_0) \dots (z_r-m_rz_0)=0.
\end{equation}
 Then the arrangement $\A(m)$ obtained as the intersection of $\hat{\A}(m)$ with 
 a generic $2$-plane $U \cong \C\p^2$ in $\C\p^r$ is a type $m$ split solvable arrangement. 
 The lines of $\A(m)$ are the traces on $U$ of the hyperplanes of $\hat{\A}(m)$. 
 
 \begin{remark}
 \label{rem:hat}
 It follows from \cite{Ha} that an arrangement $\A$ is generic if and only if
 it is the intersection of $\hat{\A}(m)$ with a generic linear subspace of $\C^{r+1}$
 of dimension at least $3$, for $m=(1,\dots, 1)$.
 \end{remark}
 
Split solvable arrangements are remarkable for illustrating a perfect 
equivalence between combinatorics and topology.

\begin{theorem}[\cite{CDP}]
\label{cor17}
For an arrangement $\A \subset \p^2$ the following are equivalent:
\begin{enumerate}
\item \label{comb}
$\A$ has the same combinatorics as a split solvable arrangement of type $m$.
\item \label{toptype}
The complement $M_{\A}$ is homeomorphic to the complement of the type $m$ arrangement $\A(m)$.                       
\end{enumerate}
\end{theorem}

This is  a consequence of a more general result concerning {\it nice} 
arrangements of lines in $\p^2$ (see \cite[Theorem 1.6]{CDP}).

\begin{remark}
\label{rem:more}
It is easy to see that $l(\A)=r$, for a split solvable arrangement of 
type $m=(m_1, \dots, m_r)$. By Theorem \ref{thm:lr}, $\A$ is not supersolvable, if $r>3$.
Clearly, $\A$ cannot be generic, if there is some $i$ such that $m_i>1$.
\end{remark}

\section{Cohomology jumping loci}
\label{sect:VR}

 Let $M$ be a connected finite type CW-complex with torsion free first homology group 
 $H_1(M):=H_1(M, \Z)$, for instance 
 the complement of a complex hyperplane arrangement. The 
 description of its {\it characteristic varieties}
\begin{equation}
\label{eq:charvar}
\VV^q_k(M, \Bbbk):= \{ \rho \in \Hom(\pi_1(M), \Bbbk^*) \;|\;\dim_{\Bbbk} H^q(M,\;_{\rho}\Bbbk) \geq k\},
\end{equation}
for an arbitrary field $\Bbbk$, is equivalent to the computation of twisted, rank one, 
$\Bbbk$-cohomology of $M$.

An analogous concept was introduced by Falk (\cite{Fa3}), under 
the name of {\em resonance varieties}:
\begin{equation}
\label{res}
\RR_{k}^{q}(M, \Bbbk) :=\{ \omega \in H^1(M, \Bbbk) \;|\; \dim_{\Bbbk}H^q(H^*(M,\Bbbk), \mu_{\omega})\geq k\},
\end{equation}
$\mu_{\omega}$ being left-multiplication by $\omega$ in the algebra $H^*(M,\Bbbk)$. 
(Here, one uses the assumption on $H_1(M)$ to check that the {\em Aomoto complex},
$(H^*(M,\Bbbk), \mu_{\omega})$, is indeed a chain complex, that is, $\omega \cdot \omega =0$.)

\begin{remark}
\label{varieties}
The sets $\VV^q_k(M, \Bbbk)$ are closed algebraic subvarieties of the algebraic torus 
$(\Bbbk^*)^n \cong \Hom(H_1(M), \Bbbk^*) \cong \Hom(\pi_1(M), \Bbbk^*)$. The 
resonance varieties $\RR_{k}^{q}(M, \Bbbk)$ are homogeneous subvarieties in $\Bbbk^n \cong H^1(M, \Bbbk)$. 

The resonance varieties in the case when $M$ is a complex algebraic variety are closely related to the classical results in algebraic geometry, starting with Castelnuovo-De Franchis Lemma and culminating with the isotropic subspace theorems due to Catanese and Bauer, see \cite{D0} for details and complete references.
\end{remark}

  We restrict our attention to the case $M=M_{\A}$, for an arbitrary arrangement $\A$,  
  and $\Bbbk=\C$, using a simplified notation
  (for example, $\RR_{k}^{q}(M):=\RR_{k}^{q}(M, \C)$).
  
 By construction, the resonance varieties are combinatorial invariants of an arrangement, 
 since $H^*(M,\Bbbk)$ depends only on the intersection lattice, see Theorem \ref{thm:OS}.
 
{ \it For characteristic varieties, combinatorial determination is a major open problem.} 
The next theorem  
provides a partial answer in this direction.

\begin{theorem}[Esnault-Schechtman-Viehweg, \cite{ESV}]
\label{germs}
The exponential map $\C^n \\ \overset{\exp}\longrightarrow (\C^*)^n$  induces 
an isomorphism of analytic germs,
\begin{equation}
\label{germs1}
(\RR_{k}^{q}(M_{\A}),0) \overset{\exp}\longrightarrow (\VV^q_k(M_{\A}),1)\, ,
\end{equation}
for all $q$ and $k$.
\end{theorem}

\begin{corollary}
\label{tangentcone}
The resonance variety $\RR_{k}^{q}(M_{\A})$ coincides with the tangent cone of 
the characteristic variety $\VV^q_k(M_{\A})$ at the point $1$,
$TC_1 (\VV^q_k(M_{\A}))$.
\end{corollary}

In \cite{Lib2}, Libgober proves that always $TC_1 (\VV^q_k(M))\subseteq \RR_{k}^{q}(M)$.
However, this inclusion may be strict, even for a smooth, quasi-projective complex
variety $M$; see \cite{DPS}.

\begin{theorem}[Arapura, \cite{A}]
\label{arapura}
Let $P$ be a smooth projective variety, $D$ a divisor of $P$, 
$M:=P\setminus D$ and assume the first Betti number of $P, \;b_1(P)$, is $0$. Then all
irreducible components of $\VV^q_k(M)$ are algebraic  subtori of
the character torus $\Hom(\pi_1(M), \C^*)$. 
\end{theorem}

(These subtori are in particular isomorphic to $(\C^*)^d$, where $d$ is the dimension
of the component.)
By the result of Arapura, the characteristic varieties of the complement $M_{\A}$
of an arrangement with $n$ hyperplanes 
(embedded in a complex projective space, by the remarks following Definition \ref{cone}) 
are unions of (possibly translated) subtori of $(\C^*)^n$. Here, we say that
an irreducible component $W$ of $\VV^q_k(M)$ is {\em non-translated} if $1 \in W$, 
and {\em translated} otherwise.

\begin{corollary}
\label{ops}
There is a bijection, induced by the exponential map, between the non-translated 
irreducible components  of $\VV^q_k(M_{\A})$ and the irreducible 
components of $\RR_{k}^{q}(M_{\A})$.
\end{corollary}

\begin{corollary}
\begin{enumerate}
\item\label{reslinear}
The resonance variety $\RR_{k}^{q}(M_\A)$ is the union of a finite number of linear subspaces in $\C^n$.
\item\label{compositiones} 
The non-translated irreducible components  of $\VV^q_k(M_{\A})$ are combinatorially
determined (see Falk-Yuzvinsky \cite{FY} for some nice combinatorial formulas).
\end{enumerate}
\end{corollary}

{\it Combinatorial determination of 
translated components is an open question, in spite of the major progress in their understanding due to \cite{D1}. Another 
open question in this direction concerns the dimensions of translated components.}

\begin{example}
\label{big}
In \cite[Example 4.1]{S0}, Suciu found examples of arrangements with translated components. 
For instance, 
the first characteristic variety $\VV^1_1(M_{\A})$ of the so called "deleted $B_3$"
arrangement, $\A \subset \C^3,\; Q(\A)=(x-z)(y-z)xyz(x-y+z)(x-y-z)(x-y)$, has
one translated component of dimension $1$.
\end{example}

Next, we focus our attention on the first resonance variety $\RR_{1}^{1}(M_{\A})$, 
that is, on the non-translated components of $\VV_{1}^{1}(M_{\A})$.
The {\it local components} of $\RR_{1}^{1}(M_{\A})$ are the irreducible components
of $\RR_{1}^{1}(M_{\A})$ that coincide with the subvarieties 
$\RR_{1}^{1}(M_{\A_{X}}) \subset \RR_{1}^{1}(M_{\A})$, for  $X\in \mathcal L_2(\A)$, 
where $\A_{X}:=\{H \in \A\;|\;X \subset H\}= \{ H_1, \dots, H_{m_X}\}$
and $m_X:=|\A_{X}| \geq 3$.  
The local components are thus indexed by the elements of rank $2$ in the lattice of
the arrangement which are intersections of at least three hyperplanes.

It is not difficult to write down the equations that define a local component.
In the notation from the beginning of Section \ref{sect:cohomology},
\begin{equation}
\label{RAx}
\RR_{1}^{1}(M_{\A_{X}})=\{ \sum_{i=1}^{m_X}x_{H_i}e_{H_i}\;|\;\sum_{i=1}^{m_X}x_{H_i}=0\}
\end{equation}
is a subspace of dimension $m_X-1$ in $ A^1(\A) \cong H^1(M_{\A})$. This 
implies at once that the dimension of a local component of $\RR_{1}^{1}(M_{\A})$ 
is at least $2$, and can be  arbitrarily large. 

On the other hand, for {\em non-local} components of $\RR_{1}^{1}(M_{\A})$ 
(i.e., components that are not local), one has the following beautiful recent result.

\begin{theorem}[Yuzvinsky, \cite{Y0}]
\label{yuzv}
The dimension of a non-local irreducible component $V$ of $\RR_{1}^{1}(M_{\A})$
satisfies the inequality $2 \leq \dim V \leq 3$.
\end{theorem}

The results below allow us to convert (combinatorial) information about resonance into 
(topological) information on characteristic varieties.

Under  "non-resonance" conditions on $z \in A^1_{\C}(\A)$, described 
by Schechtman-Terao-Varchenko in \cite{STV}, 
we have an isomorphism (see \cite{ESV}, \cite{STV}):
\begin{equation}
\label{eq:iso}
H^q(M_{\A}, \;_{\exp(z)}\C) \cong H^q(A^*_{\C}(\A), \mu_{z})
\end{equation}
In general, one has the following inequality (see Libgober-Yuzvinsky, \cite{LY}):
\begin{equation}
\label{lower}
\dim_{\C}H^q(M_{\A}, \;_{\exp(z)}\C) \geq 
\sup_{a \in 2\pi \sqrt{-1}\Z^n}\dim_{\C} H^q(A^*_{\C}(\A), \mu_{z+a})
\end{equation}

Finally, there are {\em modular inequalities}, for certain rational local systems 
(obtained by Papadima-Suciu in \cite{PS2},
extending a previous result of Cohen-Orlik from \cite{CO}):
\begin{equation}
\label{upper}
\dim_{\C}H^q(M_{\A},\; _\rho\C) \leq \dim _{\F_p}H^q(A^*_{\F_p}(\A), \mu_{z})
\end{equation}
Here $z=(m_1, \dots, m_n), \; m_i \in \Z, \;\forall i$, 
$gcd(m_1, \dots m_n)=1$, $r \in \Z_{>0}$, $p$ is a prime,
$\rho=(\exp(2 \pi \sqrt{-1}\frac{m_1}{p^r}), \dots, \exp(2 \pi \sqrt{-1}\frac{m_n}{p^r}))$,
and $\F_p$ is the prime field with $p$ elements.

Notice that both inequalities \eqref{lower} and \eqref{upper} may be strict.

\section{Milnor fiber}
\label{sect:Milnor}

In this section, we analyze Milnor fibers of homogeneous polynomials, defined in the
Introduction. When $f$ is completely reducible into distinct linear factors, 
it is the defining polynomial of an arrangement: $f=Q(\A)$. Allowing 
multiple linear factors leads to the notion of {\em multiarrangement}.

Let $\A$ be a central arrangement with homogeneous degree $n$ defining polynomial  $Q({\A})$ 
and Milnor fibration $F_{\A}\hookrightarrow M_{\A}
\stackrel{Q({\A})}{\rightarrow}\C^*$, where $F_{\A}:=Q({\A})^{-1}(1)$.
The multiplication by a primitive $n$-root of unity induces an action (of order $n$)
on the fiber $F_{\A}$, called {\it geometric monodromy}, which induces  on 
$H_*(F_{\A}, \Q)$ the {\em algebraic monodromy} action (of order $n$). We have 
a well-known equivariant decomposition 
\begin{equation}
\label{eq:ciclo}
H_q(F_{\A}, \Q)= \bigoplus_{d |n} \big( \frac{\Q[t]}{\Phi_d} \big)^{b_{q,d}(\A)},\;\forall q,
\end{equation}
where $\Phi_d$ is the $d$th cyclotomic polynomial and $b_{q,d}(\A)$ some exponents depending on $q, d$ and $\A$; see for instance
\cite{OT, L}. Equivalently, the decomposition \eqref{eq:ciclo} is encoded by
the $q$-th characteristic polynomial of the Milnor fiber, 
$\Delta_q^{\A}(t)=\Pi_{d|n}(\Phi_d)^{b_{q,d}(\A)}$.

{\it The main problem in this context is to decide whether $H_*(F_{\A}, \Q)$ is
combinatorially determined, which is open, even in degree $*=1$.} There are
interesting results by Libgober \cite{L83, LM}, but they are formulated in
non-combinatorial terms. 

In the case of line arrangements, the zeta function $Z^{\A}(t)$, which is essentially the product of $\Delta^{\A}_q(t)^{(-1)^q}$ for $q=0,1,2$, has a simple formula, see for instance Example 6.1.10 in \cite{D2}. In particular, the knowledge of $\Delta_1^{\A}(t)$ determines the characteristic polynomial of $\Delta_2^{\A}(t)$.

On the other hand, interesting information (of a combinatorial flavor) on the characteristic polynomial $\Delta_1^{\A}(t)$ may be obtained from Corollary 6.4.15 in \cite{D2}.

{\it A related question concerns the existence of torsion
in $H_*(F_{\A}, \Z)$. The answer is yes, for multiarrangements
(see Cohen-Denham-Suciu \cite{CDS}), but the problem is open for arrangements.

There is an even more ambitious goal, namely the combinatorial description of the 
algebraic monodromy.} Note that this problem is a particular case of the
combinatorial determination question for twisted cohomology, due to 
the following recurrence formula involving the exponents from \eqref{eq:ciclo}
 (see e.g. \cite{MP}):
 \begin{equation}
 \label{xyz}
 b_{q,d}(\A)+b_{q-1,d}(\A)=b_{q}(\A, \frac{1}{d})
 \end{equation}
Here $b_{q}(\A, \frac{1}{d})=\dim_{\C}H^q(M_{\A}, \, _{\rho}\C),\;\rho$ being the 
rational equimonodromical rank one system
$(\exp(2 \pi \sqrt{-1}\frac{1}{d}), \dots, \exp(2 \pi \sqrt{-1}\frac{1}{d}))$.

We will describe, following \cite{MP}, a combinatorial formula for the algebraic
monodromy action on $H_1(F_{\A}, \Q)$, for a certain class of arrangements,
defined below.

{\it The computation of the algebraic monodromy action in {\bf all} degrees
seems an extremely difficult problem, and very few examples are known.} Among them,
we mention the full combinatorial answer found by Orlik-Randell \cite{OR} for
Milnor fibers of generic arrangements (see also \cite{CS}), extended by Choudary-Dimca-Papadima in \cite{CDP}
to the class of generic sections of the arrangements $\hat{\A}(m)$ defined in
\eqref{longeq}, for an arbitrary type $m$.

It is possible to redefine the subarrangements of Coxeter arrangements of type 
$A, B$ or $D$ in a natural way by relating them to certain {\em graphs}. These are
finite graphs having at most double edges connecting two distinct vertices and 
at most one loop at each point, and edges are labeled by a $+$ or $-$ sign. 
For such a graph $\G$ on $\{ 1,\dots,l \}$ with $n$ edges, the associated arrangement 
$\A_{\G} \subset \C^{l}$ contains $n$ hyperplanes: for a $+$, 
respectively a $-$ signed edge connecting the  vertices $i \neq j$, one has a hyperplane
of equation $x_i+x_j=0$, respectively $x_i-x_j=0$, and for a loop at  $i$ 
a hyperplane of equation $x_i=0$.

We will call the arrangements obtained in this way {\it graphic} arrangements.
They coincide clearly with subarrangements of 
Coxeter arrangements of type $B$. (Note that our terminology is nonstandard:
in \cite{OT} for instance, the term "graphic arrangement" means 
"subarrangement of a Coxeter arrangement of type $A$".)
If $\G$ has no loops then it describes a subarrangement in a type $D$ Coxeter arrangement;
if we add one more restriction on $\G$, that is, to have only $"-"$ edges, we obtain 
subarrangements in a type $A$ Coxeter arrangement.

For instance, the graph $\G$ in the figure below 

%\begin{figure}
\vskip .7in
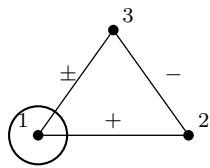
\begin{pspicture}(0,1.5)
\pscircle*(4.1,1){.07}
\pscircle(4.1,1){.4}
\pscircle*(6.1,1){.07}
\pscircle*(5.1,2.4){.07}
\psline[linewidth=.5pt](5.1,2.4)(4.1,1)(6.1,1)(5.1,2.4)
\rput(5.1,1.2){\footnotesize $+$}
\rput(4.5,1.8){\footnotesize $\pm$}
\rput(5.9,1.8){\footnotesize $-$}
\rput(3.9,1.2){\footnotesize $1$}
\rput(6.3,1.2){\footnotesize $2$}
\rput(5.3,2.6){\footnotesize $3$}
\end{pspicture}
\vskip .03in
%\end{figure}

\noindent describes the  arrangement $\A_{\G}$ in $\C^3$ of equation 
$x_1(x_1+x_3)(x_1-x_3)(x_1+x_2)(x_2-x_3)=0$.

\noindent \begin{remark}
\label{nonHyper}
Falk-Randell \cite{FR0} obtained a combinatorial formula for the Hil\-bert series of 
$\gr_*(G)\otimes \Q$, valid for fundamental groups of fiber-type arrangements. This
was extended to the hypersolvable class by Jambu-Papadima in \cite{JP1}. A new type
of combinatorial formula was found by Papadima-Suciu in \cite{PS06}, valid for the 
so-called {\em decomposable} arrangements. For a graph $\G$ with no loops or "$+$"-edges,
one knows that $\A_{\G}$ is decomposable if and only if $\G$ contains no complete subgraph
on $4$ vertices. This leads to examples of graphic arrangements which are decomposable,
but not hypersolvable; see \cite{PS06}.
\end{remark}

\begin{definition}
\label{def:baom}
For an arbitrary arrangement $\A$, the {\em Betti-Aomoto number} modulo a prime $p$ is
\[
\beta_p (\A):= \dim_{\F_p}H^1(A^*_{\F_p}(\A), \mu_{\omega_1})\, ,
\]
where $\omega_1= \sum_{H \in \A}e_H \in A^1_{\F_p}(\A)$.
\end{definition}

We may give now the promised combinatorial formula for the algebraic monodromy.

\begin{theorem}[\cite{MP}]
\label{graph1}
Let $\A$ be a graphic arrangement of rank at least $3$, with $n$
hyperplanes. Then:
\begin{equation}
\label{T1}
\Delta_1^{\A}(t)= (t-1)^{n-1}
(\Phi_2(t)\Phi_4(t))^{\beta_2(\A)} \Phi_3(t)^{\beta_3(\A)}\Phi_5(t)^{\beta_5(\A)}
\end{equation}
\end{theorem} 

We conjecture that the formula \eqref{T1} holds for subarrangements of rank at least $3$,
in arbitrary Coxeter arrangements.

Following Randell \cite{Ra}, we say that two arrangements are  {\it lattice isotopic} if
there is a continuous deformation of one arrangement into the other, through 
arrangements with the same lattice.
  
\begin{theorem}[\cite{MP}]
\label{graph2}
In Theorem \ref{graph1}, one has
\begin{center}
\label{T21}
$\Delta_1^{\A}(t)= (t-1)^{n-1}, \; when \;\A \not\equiv D_3, \;D_4$
\end{center}
or
\begin{center}
\label{T22}
$\Delta_1^{\A}(t)= (t-1)^{n-1}(t^2+t+1), \; when \;\A\equiv D_3, \; or\; D_4,$
\end{center}
where $\equiv$ means lattice isotopy.
\end{theorem}

  In the proof, the connection with  twisted cohomology is exploited, 
  together with a series of results that reduce the computation of the latter  to
  a purely combinatorial problem. The "non-resonant" case  is treated using \eqref{eq:iso},
  combined with a result of Yuzvinsky \cite{Yuz},
  concerning the cohomology of an Aomoto complex in arbitrary characteristic.
  An essential role in the resonant case is played by the modular inequalities \eqref{upper}.

 With a different approach Settepanella \cite{Set} proves that 
 the monodromy action on $H^q(F_{\A}, \Q)$ 
 is trivial for $q$ big enough, but only for full
 Coxeter arrangements of type $A,\;B$ and $D$.

\bibliographystyle{amsplain}

\end{document}